\documentclass[11pt]{article}
\usepackage{amssymb,amsthm}
\usepackage{amsmath}
\usepackage{amscd}
\usepackage{amsfonts}
\usepackage{epsf,graphics,graphicx,psfrag}

\newcommand{\des}{\displaystyle}
\newcommand{\en}{{\rm{End}}}
\newcommand{\id}{{\rm{id}}}
\newcommand{\alg}{{\rm{alg}}}
\newcommand{\sym}{{\rm{Sym}}}
\newtheorem{theorem}{Theorem}
\newtheorem{corollary}{Corollary}
\newtheorem{proposition}{Proposition} 

\newtheorem{definition}{Definition}

\begin{document}
\setlength{\baselineskip}{16pt}
\title{Symmetric quantum Weyl algebras}
\author{Rafael D\'\i az \thanks{Work partially supported by UCV.}
and Eddy Pariguan \thanks{Work partially supported by FONACIT.}}
\maketitle

\begin{abstract}
We study the symmetric powers of four algebras: $q$-oscillator
algebra, $q$-Weyl algebra, $h$-Weyl algebra and $U({\mathfrak
{sl}}_2)$. We provide explicit formulae as well as combinatorial
interpretation for the normal coordinates of products of arbitrary
elements in the above algebras.
\end{abstract}

\section{Introduction}
\vspace{0.5cm}

This paper takes part in the time-honored tradition of studying an
algebra by first choosing a ``normal'' or ``standard'' basis for
it $\mathbb{B}$, and second, writing down explicit formulae and,
if possible, combinatorial interpretation for the representation
of the product of a finite number of elements in $\mathbb{B}$, as
linear combination of elements in $\mathbb{B}$.

This method has been successfully applied to many algebras, most
prominently in the theory of symmetric functions (see \cite{GR}).
We shall deal with algebras given explicitly as the quotient of a
free algebra, generated by a set of letters $L$, by a number of
relations. We choose normal basis for our algebras by fixing an
ordering of the set of the letters $L$, and defining $\mathbb{B}$
to be the set of normally ordered monomials, i.e., monomials in
which the letters appearing in it  respect the order of $L$.

We consider algebras of the form $\sym^{n}(A)$, i.e, symmetric
powers of certain algebras. Let us recall that for each
$n\in\mathbb{N}$, there is a functor
$\sym^{n}\!:\mathbb{C}\mbox{-}
\alg\longrightarrow
\mathbb{C}\mbox{-}\alg$ from the category of associative
$\mathbb{C}$-algebras into itself defined on objects as follows:
if $A$ is a $\mathbb{C}$-algebra, then $\sym^{n}(A)$ denotes the
algebra whose underlying vector space is the $n$-th symmetric
power of $A$:\\ $\sym^{n}(A)= (A^{\otimes n})/ \langle a_1\otimes
\dots \otimes a_n -a_{\sigma^{-1} (1)}\otimes \dots \otimes
a_{\sigma^{-1}(n)}: a_i\in A, \sigma\in \mathbb{S}_n\rangle ,$
where $\mathbb{S}_n$ denotes the group of permutations on $n$
letters. The product of $m$ elements in $\sym^{n}(A)$ is given by
the rule
\begin{equation}\label{PF}
(n!)^{m-1}\prod_{i=1}^{m}\left(\overline{\bigotimes_{j=1}^{n}
a_{ij}}\right)  = \sum_{\sigma \in \{\id\} \times
\mathbb{S}_n^{m-1}}\overline{\bigotimes_{j=1}^{n}\left(
\prod_{i=1}^{m} a_{i \sigma^{-1}_i(j)}\right)}  \end{equation} for
all $(a_{ij}) \in A^{[[1,m]]\times[[1,n]]}$. Notice that if $A$ is
an algebra, then $A^{\otimes n}$ is also an algebra.
$\mathbb{S}_n$ acts on $A^{\otimes n}$ by algebra automorphisms,
and thus we have a well defined invariant subalgebra $(A^{\otimes
n})^{\mathbb{S}_n} \subseteq A^{\otimes n}$. The following result
is proven in $\cite{DP}$.
\begin{theorem} The map $s:\sym^{n}(A) \longrightarrow (A^{\otimes n})^{\mathbb{S}_n}$
given by
$$s \left(\overline{\bigotimes_{i=1}^{n} a_i}\right)
=\frac{1}{n!} \sum_{\sigma \in \mathbb{S}_n}
\overline{\bigotimes_{i=1}^{n}a_{\sigma^{-1}(i)}},\ \ \mbox{for
all}\ \ a_i\in A$$ defines an algebra isomorphism.
\end{theorem}
The main goal of this paper is the study of the symmetric powers
of certain algebras that may be regarded as quantum analogues of
the Weyl algebra. Let us recall the well-known

\begin{definition} The algebra
$\mathbb{C} \langle x,y \rangle[h]/ \langle yx-xy-h \rangle$ is
called the {\em Weyl algebra}.
\end{definition}

This algebra admits a natural representation as indicated in the
\begin{proposition}\label{rweyl}
The map $\rho:\mathbb{C} \langle x,y \rangle[h]/ \langle yx-xy-h
\rangle \longrightarrow \en(\mathbb{C}[x,h])$ given by
$\rho(x)(f)=xf$, $\rho(y)(f)=h\frac{\partial f}{\partial x}$,
$\rho(h)(f)=hf$ for any $f\in \mathbb{C}[x,h]$  defines a
representation of the Weyl algebra.
\end{proposition}

The symmetric powers of the Weyl algebra have been studied from
several point of view in papers such as
\cite{Al},\cite{EM},\cite{Etg2},\cite{Min2}. Our interest in the
subject arose from the construction of non-commutative solitonic
states in string theory, based on the combinatorics of the
annihilation $\frac{\partial}{\partial x}$ and creation $\hat{x}$
operators given in \cite{EM}. In \cite{DP} we gave explicit
formula, as well as combinatorial interpretation for the normal
coordinates of monomials $\partial^{a_1}x^{b_1}\dots\partial^{a_n}
x^{b_n}$. This formulae allow us to find explicit formulae for the
product of a finite number of elements in the symmetric powers of
the Weyl algebra. Looking at Proposition \ref{rweyl} ones notices
that the definition of Weyl algebra relies on the notion of the
derivative operator $\frac{\partial}{\partial x}$ from classical
infinitesimal calculus. The classical derivative
$\frac{\partial}{\partial x}$ admits two well-known discrete
deformations, the so called $q$-derivative $\partial_q$ and the
$h$-derivative $\partial_h$. The main topic of this paper is to
introduce the corresponding $q$ and $h$-analogues for the Weyl
algebra, and generalize the results established in \cite{DP} to
these new contexts.

For the $q$-calculus, we will actually introduce two $q$-analogues
of the Weyl algebra: the $q$-oscillator algebra (section
\ref{qos}) and the $q$-Weyl algebra (section
\ref{qweyl}). Needless to say, the $q$-oscillator algebra  also known as the $q$-boson algebra \cite{Sol},  and
$q$-Weyl algebra are deeply related. The main difference between
them is that while the $q$-oscillator is the algebra generated by
$\partial_q$ and $\hat {x}$, a third operator, the $q$-shift $s_q$
is also present in the $q$-Weyl algebra. We believe that $s_q$ is
as fundamental as $\partial_q$ and $\hat{x}$. The reason it has
passed unnoticed in the classical case is that for $q=1$, $s_q$ is
just the identity operator. For both $q$-analogues of the Weyl
algebra we are able to find explicit formulae and combinatorial
interpretation for the product rule in their symmetric powers
algebras.

For the $h$-calculus, also known as the calculus of finite
differences, we introduce the $h$-Weyl algebra in section
\ref{hweyl}. Besides the annihilator $\partial_h$ and the creator $\hat{x}$
operators, also includes an $h$-shift operator $s_h$, which again
reduces to the identity for $h=0$. We give explicit formulae and
combinatorial interpretation for the product rule in the symmetric
powers of the $h$-Weyl algebra.

In section \ref{uea} we deal with an algebra of a different sort.
Since our method has proven successful for dealing with the Weyl
algebra (and it $q$ and $h$-deformations); and it is known that
the Weyl algebra is isomorphic to the universal enveloping algebra
of the Heisenberg Lie algebra, it is an interest problem to apply
our constructions for other Lie algebras. We consider here only
the simplest case, that of ${\mathfrak {sl}}_2$. We give explicit
formulae for the product rule in the symmetric powers of
$U({\mathfrak {sl}_2})$.

Although some of the formulae in this paper are rather cumbersome,
all of them are just the algebraic embodiment of fairly elementary
combinatorial facts. The combinatorial statements will be further
analyzed in $\cite{DP1}$.

\subsection*{Notations and conventions}
\begin{itemize}
\item{$\mathbb{N}$ denotes the set of natural numbers. For $x\in
\mathbb{N}^{n}$ and $i\in\mathbb{N}$, we denote by $x_{<i}$ the
vector $(x_1,\dots,x_{i-1})\in \mathbb{N}^{i-1}$, by $x_{\leq i}$
the vector $(x_1,\dots,x_{i})\in \mathbb{N}^{i}$, by $x_{>i}$ the
vector $(x_{i+1},\dots,x_{n})\in \mathbb{N}^{n-i}$ and by $x_{\geq
i}$ the vector $(x_i,\dots,x_n)\in \mathbb{N}^{n-i+1}$.}
 \item
{Given $(U,<)$ an ordered set and $s\in U$, we set $U_{>
s}:=\{u\in U: u>s\}$.}
 \item{ $|\mbox{
}|:\mathbb{N}^{n}\longrightarrow \mathbb{N}$ denotes the map such
that $|x|:=\sum_{i=1}^{n} x_i$, for all $x\in \mathbb{N}^{n}$.}
\item{For a set $X$, $\sharp(X):=$ cardinality of $X$, and
$\mathbb{C}\langle X \rangle\!\!:=$ free associative algebra
generated by $X$.} \item{Given natural numbers $a_1,\dots,a_n\in
\mathbb{N}$, $\min(a_1,\dots,a_n)$ denotes the smallest number in
the set $\{a_1,\dots,a_n\}$.}
\item{Given $n\in\mathbb{N}$, we set $[[1,n]]=\{1,\dots,n\}$.}

 \item{Let $S$ be a set and
$A:[[1,m]]\times [[1,n]]\longrightarrow S$ an $S$-valued matrix.
For $\sigma \in (\mathbb{S}_n)^{m}$ and $j\in [[1,n]]$,
$A_j^{\sigma}:[[1,m]]\longrightarrow S$ denotes the map such that
$A_{j}^{\sigma}(i)=A_{i \sigma^{-1}_i(j)}$, for all
$i=1,\dots,m.$} \item{The $q$-analogue $n$ integer is ${\des
[n]:=\frac{1-q^{n}}{1-q}}$. For $k\in \mathbb{N}$, we will use
$[n]_{k}:=[n][n-1]\dots [n-k+1]$.}
\end{itemize}

\section{Symmetric $q$-oscillator algebra}\label{qos}
\vspace{0.5cm}
In this section we define the $q$-oscillator algebra and study its
symmetric powers. Let us introduce several fundamental operators
in $q$-calculus (see
$\cite{Kac}$ for a nice introduction to $q$-calculus).\\

\begin{definition}
The operators $\partial_q, s_q, \hat{x}, \hat{q}, \hat{h}:
\mathbb{C}[x,q,h] \longrightarrow \mathbb{C}[x,q,h]$ are given as
follows
 $$\begin{array}{lclcc}
  \partial_q f(x)=\frac{f(qx)-f(x)}{(q-1)x} &\mbox{} & s_q(f)(x)=f(qx)
  & \mbox{ } &\hat{q}(f)=qf\\
  \hat{x}(f)=xf &\mbox{ } & \hat{h}(f)=hf &\mbox{ } &  \\
 \end{array}$$ for all $f\in \mathbb{C}[x,q,h]$. We call $\partial_q$
 the $q$-derivative and $s_q$ the $q$-shift.
\end{definition}

\begin{definition}
The algebra $\mathbb{C}\langle x,y\rangle[q,h]/I_{qo}$, where
$I_{qo}$ is the ideal generated by the relation $yx=qxy+h$ is
called the {\em $q$-oscillator algebra}.
\end{definition}
We have the following analogue of Proposition $\ref{rweyl}$.
\begin{proposition}\label{rqo}
The map $\rho:\mathbb{C}\langle x,y\rangle[q,h]/I_{qo}
\longrightarrow \en(\mathbb{C}[x,q,h])$ given by
$\rho(x)=\hat{x}$, $\rho(y)=h\partial_q$, $\rho(q)=\hat{q}$ and
$\rho(h)=\hat{h}$ defines a representation of the $q$-oscillator
algebra.
\end{proposition}
Notice that if we let $q\rightarrow 1$, $y$ becomes central and we
recover the Weyl algebra. We order the letters of the
$q$-oscillator algebra as follows:
$q<x<y<h$.\\

Assume we are given $A=(A_1,\dots,A_n)\in (\mathbb{N}^{2})^{n}$,
and $A_i=(a_i,b_i)\in \mathbb{N}^{2}$, for $i\in [[1,n]]$. Set
$X^{A_{i}}=x^{a_i}y^{b_i}$, for $i\in [[1,n]]$. We set
$a=(a_1,\dots, a_n)\in \mathbb{N}^{n}$, $b=(b_1,\dots, b_n)\in
\mathbb{N}^{n}$, $c=(c_1,\dots, c_n)\in \mathbb{N}^{n}$ and
$|A|=(|a|,|b|)\in \mathbb{N}^{2}$. Using
this notation we have\\

\begin{figure}[h]
\begin{center}
\psfrag{t}{${\scriptscriptstyle{\infty}}$}
\psfrag{a}{${\scriptscriptstyle {U_1}}$}
\psfrag{b}{${\scriptscriptstyle {V_1}}$}
\psfrag{c}{${\scriptscriptstyle {U_2}}$}
\psfrag{d}{${\scriptscriptstyle {V_2}}$}
\psfrag{f}{${\scriptscriptstyle {U_3}}$}
\psfrag{g}{${\scriptscriptstyle{V_3}}$}
\includegraphics[height=1.3cm]{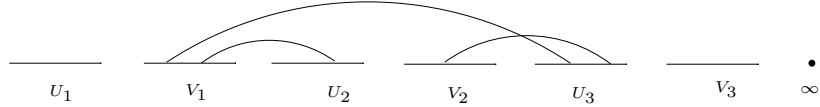}
\caption{Combinatorial interpretation of
$N_{qo}$.}\label{fg:crnum}
\end{center}
\end{figure}

\begin{definition}
The {\em normal coordinates} $ N_{qo}(A,k)$ of
${\des\prod_{i=1}^{n}X^{A_{i}}}\in \mathbb{C}\langle
x,y\rangle[h,q]/I_{qo}$ are given by the identity
\begin{equation}\label{cqo}
{\des\prod_{i=1}^{n}X^{A_{i}}}=
 {\des \sum_{k=0}^{\min} N_{qo}(A,k) X^{|A|-(k,k)} h^{k}}
\end{equation}
where $\min=\min(|a|,|b|)$. For $k>\min$, we set $N_{qo}(A,k)$ to
be equal to $0$.
\end{definition}
Let us introduce some notation needed to formulate Theorem
$\ref{ncqo}$ below which provides explicit formula for the normal
coordinates  $N_{qo}(A,k)$ of ${\des\prod_{i=1}^{n}X^{A_{i}}}$.
Given $(A_1,\dots,A_n)\in (\mathbb{N}^{2})^{n}$ choose disjoint
totally ordered sets $(U_i,<_{i})$, $(V_i,<_{i})$ such that
$\sharp(U_i)=a_i$ and $\sharp (V_i)=b_i$, for $i\in [[1,n]]$.
Define a total order set $(U\cup V\cup\{ \infty\},<)$, where
${\des U=\cup_{i=1}^{n} U_i}$, $V={\des\cup_{i=1}^{n}V_i}$ and
$\infty\not\in U\cup V$, as follows: Given $u,v\in U\cup V\cup\{
\infty\}$ we say that $u\leq v$ if and only if a least one of the
following conditions hold
\begin{eqnarray*}
 u\in V_i, v\in V_j\ \mbox{and}\ i\leq j; &\mbox{ }&
u\in V_i, v\in U_j\ \mbox{and}\ i\leq j;\\
 u\in U_i, v\in V_j \ \mbox{and}\ i\leq j;  &\mbox{}&
u,v\in U_i \ \mbox{and}\  u\leq_i v;\\
 u,v\in V_i \ \mbox{and}\  u\leq_i v;  \mbox{ } \mbox{ } \mbox{ } \mbox{ }
\mbox{ }  &\mbox{ }& v=\infty.
\end{eqnarray*}
Given $k\in \mathbb{N}$, we let $P_{k}(U,V)$ be the set of all
maps $p:V\longrightarrow U\cup\{\infty\}$ such that
\begin{itemize}
\item{$p$ restricted to $p^{-1}(U)$ is injective, and $\sharp
(p^{-1}(U))=k$.} \item{If $(v,p(v))\in V_i\times U_j$ then $i<j$.}
\end{itemize}
Figure $\ref{fg:crnum}$ shows an example of such a map. We only
show the finite part of $p$, all other points in $V$ being mapped
to $\infty$.

\begin{definition}
The value of the {\em crossing number map}
$c:P_{k}(U,V)\longrightarrow \mathbb{N}$ when evaluated on $p\in
P_k(U,V)$ is given by \vspace{-0.1cm}
$$c(p)=\sharp (\{(s,t)\in V\times U|  s<t<p(s), t\in p(U_{> s})^{c}\}).$$
\end{definition}
\begin{theorem}\label{ncqo}
For any $A=(A_1,\dots,A_n)\in (\mathbb{N}^{2})^{n}$ with
$A_i=(a_i,b_i)\in\mathbb{N}^{2}$ for $i\in [[1,n]]$ and any $k\in
\mathbb{N}$, we have that
$$N_{qo}(A,k)={\des\sum_{p\in P_k(U,V)} q^{c(p)}}.$$
\end{theorem}
\begin{proof}
The proof is by induction. The only non-trivial case is the
following
\begin{eqnarray}\label{recur} {\des y \prod_{i=1}^{n}
X^{A_i}}& =& {\des \sum_{k=0}^{\min} q^{|a|-k}
N_{qo}(A,k) X^{|A|-(k,k+1)} h^{k}} \nonumber \\
\mbox{} & \mbox{} & +{\des \sum_{k=0}^{\min} N_{qo}(A,k)
\sum_{i=1}^{a-k} q^{i-1} X^{|A|-(k+1,k)}h^{k}}
\end{eqnarray}
where $\min=\min(|a|,|b|)$. Normalizing the left-hand side of
$(\ref{recur})$ we get a recursive relation
$$N_{qo}(((0,1)A),k)=N_{qo}(A,k)+\sum_{i=1}^{a-k}q^{i-1} N_{qo}(A,k-1). $$
The other recursion needed being $N_{qo}(((1,0)A),k)=N_{qo}(A,k)$.
It is straightforward to check that $N_{qo}(A,k)={\des\sum_{p\in
P_k(U,V)} q^{c(p)}}$ satisfies both recursions.
\end{proof}
\vspace{0.3cm} Consider the identity $(\ref{cqo})$ in the
representation of the $q$-oscillator algebra defined in
Proposition $\ref{rqo}$. Apply both sides of the identity
$(\ref{cqo})$ to $x^{t}$ for $t\in \mathbb{N}$, and use Theorem
$\ref{ncqo}$ to get the fundamental \vspace{0.3cm}
\begin{corollary} Given $(t,a,b)\in \mathbb{N}\!\times\! \mathbb{N}^{n}\!\times\! \mathbb{N}^{n}$ the following
identity holds $${\des
\prod_{i=1}^{n}[t+|a_{>i}|-|b_{>i}|]_{b_i}=\sum_{p\in P_k(U,V)}
q^{c(p)}[t]_{|b|-k} }.$$
\end{corollary}
Our next theorem gives a fairly simple formula for the product of
$m$ elements in $\sym^{n}(\mathbb{C}\langle
x,y\rangle[q,h]/I_{qo})$. Fix a matrix $A:[[1,m]]\!\times\!
[[1,n]]\longrightarrow \mathbb{N}^{2}$,
$(A_{ij})=((a_{ij}),(b_{ij}))$. Recall that given $\sigma \in
(\mathbb{S}_n)^{m}$ and $j\in[[1,n]]$, $A_j^{\sigma}$ denotes the
vector $(A_{1
\sigma_{1}^{-1}(j)},\dots,A_{m\sigma_{m}^{-1}(j)})\in
(\mathbb{N}^{2})^{m}$ and $X_j^{A_{ij}}=x_j^{a_{ij}}
y_j^{b_{ij}}$. Set \\ $|A_j^{\sigma}|=(|a_j^{\sigma}|,
|b_j^{\sigma}|)$ where ${\des |a_j^{\sigma}|=\sum_{i=1}^{m}
a_{i\sigma^{-1}_i(j)}}$ and ${\des |b_j^{\sigma}|=\sum_{i=1}^{m}
b_{i\sigma^{-1}_i(j)}}$.
\begin{theorem} For any $A:[[1,m]]\times [[1,n]]\longrightarrow
\mathbb{N}^{2}$, the identity
$${\des (n!)^{m-1}\prod_{i=1}^{m}  \overline{\prod_{j=1}^{n}X_j^{A_{ij}}}=
\sum_{\sigma,k}\left( \prod_{j=1}^{n} N_{qo}(A_j^{\sigma}
,k_j)\right)\overline{\prod_{j=1}^{n}
X_j^{|A_j^{\sigma}|-(k_j,k_j)}}h^{|k|} }$$ where $\sigma\in \{
\id\} \times \mathbb{S}_n^{m-1}$ and $k\in\mathbb{N}^{n}$, holds
in $\sym^{n}(\mathbb{C}\langle x,y\rangle[q,h]/I_{qo})$.
\end{theorem}
\begin{proof} Using the product rule given in $(\ref{PF})$, the identity
$(\ref{cqo})$ and the distributive property we obtain
\begin{eqnarray*}
(n!)^{m-1}\prod_{i=1}^{m} \overline{\prod_{j=1}^{n}
X_j^{A_{ij}}}&=&\sum_{\sigma \in {\{\id\}
\times \mathbb{S}_n^{m-1}}}\prod_{j=1}^{n} \overline{\prod_{i=1}^{m} X_j^{A_{i \sigma_{i}^{-1}(j)}}}\\
&=&\sum_{\sigma \in {\{\id\} \times
\mathbb{S}_n^{m-1}}}\overline{\prod_{j=1}^{n} \left( \sum_{k=0}^{\min_j}
N_{qo}(A_j^{\sigma}
,k) x_j^{|a_j^{\sigma}|-k} y_j^{|b_j^{\sigma}|-k} h^{k}\right)}\\
&=& \sum_{\sigma,k}\left( \prod_{j=1}^{n} N_{qo}(A_j^{\sigma}
,k_j)\right)\overline{\prod_{j=1}^{n}
X_j^{|A_j^{\sigma}|-(k_j,k_j)}}h^{|k|}
\end{eqnarray*}
where $\min_j=\min(|a_j^{\sigma}|,|b_j^{\sigma}|)$.
\end{proof}

\section{Symmetric $q$-Weyl algebra}\label{qweyl}

\vspace{0.5cm} In this section we study the symmetric powers of
the $q$-Weyl
algebra.\\

\begin{definition}
The $q$-Weyl algebra is given by $\mathbb{C}\langle x,y,z
\rangle[q]/I_q$, where $I_q$ is the ideal generated by the
following relations:

$$\begin{array}{ccccc} zx=xz+y &\mbox{ } & yx=qxy & \mbox{}  &zy=qyz \\
\end{array}$$
\end{definition}
 We have the following $q$-analogue of Proposition
$\ref{rweyl}$.
 \begin{proposition}\label{rpq}
The map $\rho:\mathbb{C}\langle x,y,z \rangle[q]/I_q
\longrightarrow \en(\mathbb{C}[x,q])$ given by $\rho(x)=\hat{x}$,
$\rho(y)=s_q$, $\rho(z)=\partial_q$ and $\rho(q)=\hat{q}$ defines
a representation of the $q$-Weyl algebra.
 \end{proposition}
Notice that if we let $q\rightarrow 1$, we recover the Weyl
algebra. We order the letters of the $q$-Weyl algebra as follows:
$q<x<y<z$. Given $a\in \mathbb{N}$ and $I\subset [[1,a]]$, we
define the {\em crossing number of $I$} to be
$$\chi(I):=\sharp(\{(i,j): i>j, i\in I, j\in I^{c}\}).$$ For $k\in \mathbb{N}$, we let
$\chi_k:\mathbb{N}\longrightarrow \mathbb{N}$  the map given by
${\des \chi_k(a)=\sum_{\begin{array}{c}
  {\scriptstyle {\sharp (I)=k}} \\
  {\scriptstyle {I\subset [[1,a]]}}\\
\end{array}} q^{\chi(I)}}$, for all $a\in \mathbb{N}$. We have the following
\vspace{0.3cm}
\begin{theorem}\label{qfor} Given $a,b\in \mathbb{N}$, the
following identities hold in $\mathbb{C}\langle x,y,z
\rangle[q]/I_q $
\begin{enumerate}
\item{$z^{a}x^{b}={\des \sum_{k=0}^{\min}\chi_k(a)[b]_{k} x^{b-k}
y^{k} z^{a-k} }$, where $\min=\min(a,b)$.}
\item{$z^{a}y^{b}=q^{ab}y^{b}z^{a}$.}
\item{$y^{a}x^{b}=q^{ab}x^{b}y^{a}$.}
\end{enumerate}
\end{theorem}
\begin{proof}
$2.$ and $3.$ are obvious, let us to prove $1.$ It should be clear
that
$${\des{z^{a}x^{b}=\sum_{I\subset
[[1,a]]}[b]_{\sharp(I)}x^{b-\sharp(I)}\prod_{j=1}^{a} f_I(j)}},$$
where $f_I(j)=z$, if $j\notin I$ and $f_I(j)=y$, if $j\in I$. The
normal form of ${\des\prod_{j=1}^{a} f_I(j)}$ is $q^{\chi(I)}
y^{\sharp(I)} z^{a-\sharp(I)}$. Thus,\\
$$z^{a}x^{b}={\des\sum_{k=0}^{\min}\left( \sum_{I}
q^{\chi(I)}\right) [b]_k x^{b-k}y^{k}z^{a-k}
=\sum_{k=0}^{\min}\chi_k(a)[b]_k x^{b-k}y^{k}z^{a-k}},$$ where
$\min=\min(a,b)$,  $I\subset [[1,a]]$ and $\sharp (I)=k$.
\end{proof}
\vspace{0.3cm}
Assume we are given $A=(A_1,\dots,A_n)\in (\mathbb{N}^{3})^{n}$,
where $A_i=(a_i,b_i,c_i)$, for $i\in[[1,n]]$, also
$X^{A_{i}}=x^{a_i}y^{b_i}z^{c_i}$ for $i\in[[1,n]]$. We set
$a=(a_1,\dots,a_n)\in \mathbb{N}^{n}$, $b=(b_1,\dots,b_n)\in
\mathbb{N}^{n}$, $c=(c_1,\dots,c_n)\in \mathbb{N}^{n}$ and
$|A|=(|a|,|b|,|c|)\in \mathbb{N}^{3}$. Using this notation, we
have the\\
\begin{definition}
The {\em normal coordinates} $ N_{q}(A,k)$ of
${\des\prod_{i=1}^{n}X^{A_i}}\in \mathbb{C}\langle x,y,z
\rangle[q]/I_q$ are given via the identity
\begin{equation}\label{cq}
{\des\prod_{i=1}^{n}X^{A_i}= \sum_{k\in \mathbb{N}^{{n-1}} }
N_{q}(A,k) X^{|A|+(-|k|,|k|,-|k|)}}
\end{equation}
where $k$ runs over all vectors $k=(k_1,\dots,k_{n-1})\in
\mathbb{N}^{n-1}$ such
that $0\leq k_i \leq \min (|c_{\leq i}|-|k_{<i}|,a_{i+1})$. We set
$N_{q}(A,k)=0$ for $k\in\mathbb{N}^{n-1}$ not satisfying the
previous conditions.
\end{definition}
Our next theorem follows from Theorem $\ref{qfor}$ by induction.
It gives an explicit formula for the normal coordinates
$N_{q}(A,k)$ in the $q$-Weyl algebra of
${\des\prod_{i=1}^{n}X^{A_i}}$.
\begin{theorem}\label{nq} Let $A,k$  be as in the previous
definition, we have
$$N_{q}(A,k)={\des q^{\sum_{i=1}^{n-1} \lambda(i)}\prod_{j=1}^{n-1}
\chi_{k_j}(|c_{\leq j}|-|k_{<j}|) [a_{j+1}]_{k_{j}}},$$ where
$\lambda(i)=b_{i+1}(|c_{\leq i}|-|k_{\leq
i}|)+(a_{i+1}-k_i)(|b_{\leq i}|+|k_{<i}|)$.
\end{theorem}
Applying both sides of the identity $(\ref{cq})$ in the
representation of the $q$-Weyl algebra given in Proposition
$\ref{rpq}$ to $x^{t}$ and using Theorem $\ref{nq}$, we obtain the
remarkable
\begin{corollary} For any given
$(t,a,b,c)\in\mathbb{N}\!\times\!\mathbb{N}^{n}\!\times\! \mathbb{N}^{n}\!\times\!\mathbb{N}^{n}$,
the following identity holds
$${\des\prod_{i=1}^{n}q^{\gamma(i)} [t+|a_{>i}|+|c_{>i}|]_{c_i}}
={\des \sum_{k}q^{\beta(k)}\left(
\prod_{j=1}^{n-1}\chi_{k_j}(|c_{\leq
j}|-|k_{<j}|)[a_{j+1}]_{k_j}\right)[t]_{|c|-|k|}}$$ where $k\in
\mathbb{N}^{n-1}$ such that $0\leq k_i \leq \min (|c_{\leq
i}|-|k_{<i}|,a_{i+1})$,\\
 $\gamma(i)=b_i(t+|a_{>i}|-|c_{>i-1}|)$, and
${\des\beta(k)=\left(\sum_{i=1}^{n-1}\lambda(i)\right)+(|b|+|k|)(t-|c|+|k|)}$.

\end{corollary}
Our next theorem provides an explicit formula for the product of
$m$ elements in the algebra $\sym^{n}(\mathbb{C}\langle x,y,z
\rangle[q]/I_q)$. Fix $A:[[1,m]]\times
[[1,n]]\longrightarrow \mathbb{N}^{3}$, with
$(A_{ij})=((a_{ij}),(b_{ij}),(c_{ij}))$. Recall that given
$\sigma\in \mathbb{S}_n^{m}$ and $j\in[[1,n]]$, $A_j^{\sigma}$
denotes the vector $(A_{1\sigma_1^{-1}(j)},\dots,
A_{m\sigma_m^{-1}(j)})\in (\mathbb{N}^{3})^{m}$ and set
$X_j^{A_{ij}}=x_j^{a_{ij}} y_j^{b_{ij}} z_j^{c_{ij}}$, for
$j\in[[1,n]]$. Set $A_j^{\sigma}=(|a_j^{\sigma}|,
|b_j^{\sigma}|,|c_j^{\sigma}|)$, where ${\des
|a_j^{\sigma}|=\sum_{i=1}^{m} a_{i\sigma^{-1}_i(j)}}$ and
similarly for $|b_j^{\sigma}|$ and $|c_j^{\sigma}|$. We have the
following:

\begin{theorem}\label{jo} For any $A:[[1,m]]\times [[1,n]]\longrightarrow
\mathbb{N}^{3}$, the identity
$${\des (n!)^{m-1}\prod_{i=1}^{m} \overline{ \prod_{j=1}^{n} X_j^{A_{ij}}}=
\sum_{\sigma,k}\left( \prod_{j=1}^{n}
N_{q}(A_j^{\sigma},k^{j})\right)\overline{ \prod_{j=1}^{n}
X_j^{|A_j^{\sigma}|+(-|k^{j}|,|k^{j}|,-|k^{j}|)}}}$$ where
$\sigma\in\{\id\} \times \mathbb{S}_n^{m-1}$ and
$k=(k^{1},\dots,k^{n})\in(\mathbb{N}^{m-1})^{n}$, holds in\\
$\sym^{n}(\mathbb{C}\langle x,y,z \rangle[q]/I_q)$.
\end{theorem}
\begin{proof}
\begin{eqnarray*}
(n!)^{m-1}\prod_{i=1}^{m} \overline{ \prod_{j=1}^{n}
X_j^{A_{ij}}}&=&\sum_{\sigma \in {\{\id\}
\times \mathbb{S}_n^{m-1}}}\prod_{j=1}^{n}\overline{ \prod_{i=1}^{m} X_j^{A_{i \sigma_{i}^{-1}(j)}}}\\
&=&\sum_{\sigma}\overline{\prod_{j=1}^{n} \left(
\sum_{k=0}^{\min_j}
N_{q}(A_j^{\sigma},k) x_j^{|a_j^{\sigma}|-|k|} y_j^{|b_j^{\sigma}|+|k|} z_j^{|c_j^{\sigma}|-|k|} \right)}\\
&=&\sum_{\sigma,k}\left( \prod_{j=1}^{n}
N_{q}(A_j^{\sigma},k^{j})\right)\overline{ \prod_{j=1}^{n}
X_j^{|A_j^{\sigma}|+(-|k^{j}|,|k^{j}|,-|k^{j}|)}},
\end{eqnarray*}
where $\min_j=\min(|a_j^{\sigma}|,|b_j^{\sigma}|,|c_j^{\sigma}|)$.
\end{proof}

\section{Symmetric $h$-Weyl algebra}\label{hweyl}
\vspace{0.5cm}
In this section we introduce the $h$-analogue of the Weyl algebra
in the $h$-calculus, and study its symmetric powers. A basic introduction to $h$-calculus
may be found in $\cite{Kac}$.\\

\begin{definition}
The operators $\partial_h, s_h, \hat{x}, \hat{h}: \mathbb{C}[x,h]
\longrightarrow \mathbb{C}[x,h]$ are given by \vspace{-0.2cm}
 $$\begin{array}{ccccccc}
   \partial_h
 f(x)=\frac{f(x+h)-f(x)}{h} & \mbox{ } & s_h(f)(x)=f(x+h) & \mbox{ }& \hat{x}(f)=xf & \mbox{ }  &\hat{h}(f)=hf \\
 \end{array}$$for $f\in \mathbb{C}[x,h]$. We call $\partial_h$
 the $h$-derivative and $s_h$ the $h$-shift.
\end{definition}

\begin{definition}
The {\em $h$-Weyl algebra} is the algebra $\mathbb{C}\langle x,y,z
\rangle[h]/I_h$, where $I_h$ is the ideal generated by the
following relations:
$$\begin{array}{cccccc}
 yx=xy+z   & \mbox{ } &   zx=xz+zh &  \mbox{ } & yz=zy \\
\end{array}$$
\vspace{-0.5cm}
 \end{definition}
 \begin{proposition}
The map $\rho:\mathbb{C}\langle x,y,z \rangle[h]/I_h
\longrightarrow \en(\mathbb{C}[x,h])$ given by $\rho(x)=\hat{x}$,
$\rho(y)=\partial_h$, $\rho(z)=s_h$ and $\rho(h)=\hat{h}$ defines
a representation of the $h$-Weyl algebra.
 \end{proposition}
Notice that if we let $h\rightarrow 0$, $z$ becomes a central
element and we recover the Weyl algebra. We order the letters on
the $h$-Weyl algebra as follows
 $x<y<z<h$. Also, for $a\in\mathbb{N}$ and
$k=(k_1,\dots, k_n)\in \mathbb{N}^{n}$,  we set ${\des {a\choose
k}:=\frac{a!}{\prod k_i ! (a-|k|)!}}$. With this notation, we have
the
\vspace{0.5cm}
\begin{theorem}\label{hwe} Given $a,b\in \mathbb{N}$, the
following identities holds in $\mathbb{C}\langle x,y,z
\rangle[h]/I_h$
\begin{enumerate}
\item{$z^{a}x^{b}={\des \sum_{k=0}^{b} {b \choose k} a^{k} x^{b-k} z^{a} h^{k}}$.}
\item{$z^{a}y^{b}=y^{b} z^{a}$.}
\item{$y^{a}x^{b}={\des \sum_{\scriptstyle{k\in \mathbb{N}^{a}}}
{b \choose k} x^{b-|k|} y^{a-s(k)}z^{s(k)}h^{|k|-s(k)}}$,  where
$k\in \mathbb{N}^{a}$ is such that \\ $0\leq |k| \leq b$ and
$s(k)=\sharp (\{i: k_i\neq 0\})$.}
\end{enumerate}
\end{theorem}
\begin{proof}
$2.$ is obvious and $3.$ is similar to $1.$ We prove $1.$ by
induction. It is easy to check that $z^{a}x=xz^{a}+az^{a}h$.
Furthermore,
\begin{eqnarray*}
z^{a}x^{b+1}&=&{\des \sum_{k=0}^{b} {b \choose k} a^{k} x^{b-k} (z^{a}x) h^{k}}\\
\mbox{}&=&{\des \sum_{k=0}^{b} {b \choose k} a^{k} x^{b+1-k} z^{a} h^{k}+
\sum_{k=1}^{b+1} {b \choose k-1} a^{k} x^{b+1-k} z^{a} h^{k}} \\
\mbox{}&=&  {\des \sum_{k=0}^{b+1} {b+1 \choose k} a^{k} x^{b+1-k} z^{a} h^{k}}.
\end{eqnarray*}
\end{proof}

Assume we are given $A=(A_1,\dots,A_n)\in (\mathbb{N}^{3})^{n}$
and $A_i=(a_i,b_i,c_i)$, for $i\in[[1,n]]$. Set
$X^{A_{i}}=x^{a_i}y^{b_i}z^{c_i}$ for $i\in[[1,n]]$. We set
$a=(a_1,\dots,a_n)\in \mathbb{N}^{n}$, $b=(b_1,\dots,b_n)\in
\mathbb{N}^{n}$, $c=(c_1,\dots,c_n)\in \mathbb{N}^{n}$ and
$|A|=(|a|,|b|,|c|)\in \mathbb{N}^{3}$. Using this notation, we
have the\\

\begin{definition}
The {\em normal coordinates} $ N_{h}(A,p,q)$ of
${\des\prod_{i=1}^{n}X^{A_i}}\in \mathbb{C}\langle x,y,z
\rangle[h]/I_h$ are given via the identity
\begin{equation}\label{cnh}
{\des\prod_{i=1}^{n}X^{A_i}}=
 {\des \sum_{p,q} N_{h}(A,p,q) X^{r(A,p,q)}h^{|q|+|p|-s(p)} }
\end{equation}
where the sum runs over all vectors
$q=(q_{1},\dots,q_{{n-1}})\in\mathbb{N}^{n-1}$ such that\\  $0\leq
q_{j} \leq a_{j+1}$  and $p=(p_{1},\dots,p_{{n-1}})$ with
$p_{j}\in
\mathbb{N}^{|b_{\leq j}|-|s(p_{<j})|}$. \\ Also,
$r(A,p,q)=(|a|-|p|-|q|,|b|-s(p),|c|+s(p))$, where
$s(p)=\sum_{j=1}^{n-1} s(p_{j})$. We set $N_{h}(A,p,q)=0$ for
 $p,q$ not satisfying the previous conditions.
\end{definition}
The condition for $p$ and $q$ in the definition above might seem
unmotivated. They appear naturally in the course of the proof of
Theorem \ref{ns} below, which is proved using induction and
Theorem
\ref{hwe}.

\begin{figure}[ht]
\begin{center}
\epsfxsize 11cm \epsfbox{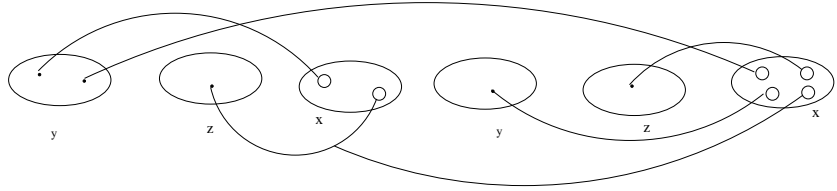}
\caption{Combinatorial interpretation of
$N_{h}$.}\label{fg:crossing}
\end{center}
\end{figure}

\begin{theorem}\label{ns} Let $A,p$ and $q$ be as in the previous
definition, we have $$N_h(A,p,q)={\des \prod_{i=1}^{n-1} {a_{i+1}
\choose q_{i}} {a_{i+1}-q_{i} \choose p_{i}}(|c_{\leq i}|+|s(p_{<i})|)^{q_i} }.$$
\end{theorem}
Figure $\ref{fg:crossing}$ illustrates the combinatorial
interpretation of Theorem $\ref{ns}$. As we try to moves the $z$'s
or the $y$'s above the `$x$', a subset of the $x$ may get killed.
The $z$'s do not die in this process but the $y$'s do turning
themselves into $z$' s.

Our next result provides an explicit formula for the product of
$m$ elements in the algebra $\sym^{n}(\mathbb{C}\langle x,y,z
\rangle[h]/I_h)$. Fix $A:[[1,m]]\times
[[1,n]]\longrightarrow \mathbb{N}^{3}$, with
$(A_{ij})=((a_{ij}),(b_{ij}),(c_{ij}))$. Recall that given
$\sigma\in \mathbb{S}_n^{m}$ and $j\in[[1,n]]$, $A_j^{\sigma}$
denotes the vector $(A_{1\sigma_1^{-1}(j)},\dots,
A_{m\sigma_m^{-1}(j)})\in (\mathbb{N}^{3})^{m}$ and set
$X_j^{A_{ij}}=x_j^{a_{ij}} y_j^{b_{ij}} z_j^{c_{ij}}$, for
$j\in[[1,n]]$. Set $A_j^{\sigma}=(|a_j^{\sigma}|,
|b_j^{\sigma}|,|c_j^{\sigma}|)$, where ${\des
|a_j^{\sigma}|=\sum_{i=1}^{m} a_{i\sigma^{-1}_i(j)}}$ and
similarly for $|b_j^{\sigma}|$ and $|c_j^{\sigma}|$.

\begin{theorem}\label{fph} For any $A:[[1,m]]\times
[[1,n]]\longrightarrow \mathbb{N}^{3}$, the identity
$${\des (n!)^{m-1}\prod_{i=1}^{m} \overline{ \prod_{j=1}^{n}
X_j^{A_{ij}}}= \sum_{\sigma,p,q}\left( \prod_{j=1}^{n}
N_{h}(A_j^{\sigma},p^{j},q^{j})\right) \overline{ \prod_{j=1}^{n}
X_j^{r(A_j^{\sigma},p^{j},q^{j})}} h^{k_j(p,q)}}$$ holds in
$\sym^{n}(\mathbb{C}\langle x,y,z \rangle/I_h)$, where $\sigma\in
\{\id\}\times \mathbb{S}_n^{m-1}$, $p^{j}$, $q^{j}$ are such that
$(A_j^{\sigma},p^{j},q^{j})$ satisfy the condition of the
definition above.
$r(A_j^{\sigma},p^{j},q^{j})=(|a_j^{\sigma}|-|p^{j}|-|q^{j}|,|b_j^{\sigma}|-s(p^{j}),|c_j^{\sigma}|+s(p^{j}))$
and $k_j(p,q)=|q^{j}|-|p^{j}|- s(p^{j})$.
\end{theorem}
Theorem \ref{fph} is proven similarly to Theorem $\ref{jo}$.

\section{Deformation quantization of $({\mathfrak {sl}}_2^{\ast})^{n}/{\scriptstyle
{\mathbb{S}_n}}$}\label{uea}
\vspace{0.5cm}

We denote by ${\mathfrak {sl}}_2$ the Lie algebra of all $2\times
2$ complex matrices of trace zero. ${\mathfrak {sl}}_2^{\ast}$ is
the dual vector space. It carries a natural structure of Poisson
manifold. We consider a deformation quantization of the Poisson
orbifold $({\mathfrak {sl}}_2^{\ast})^{n}/\mathbb{S}_n$. It is
proven in
\cite{Kon} that the quantized algebra of the Poisson manifold
${\mathfrak {sl}}_2^{\ast}$ is isomorphic to $U({\mathfrak
{sl}}_2)$ the universal enveloping algebra of ${\mathfrak
{sl}}_2$, after setting the formal parameter $\hbar$ appearing in
\cite{Kon} to be $1$. Thus we regard $(U({\mathfrak {sl}}_2)^{\otimes
n})^{\mathbb{S}_n}
\cong \sym^{n}( U({\mathfrak {sl}}_2))$ as the quantized algebra
associated to the Poisson orbifold $({\mathfrak
{sl}}_2^{\ast})^{n}/\mathbb{S}_n$. It is well-known that
$U({\mathfrak {sl}}_2)$  can be identified with the algebra
$\mathbb{C}\langle x,y,z \rangle/I_{{\mathfrak {sl}}_2}$ where
$I_{{\mathfrak {sl}}_2}$ is the ideal generated by the following
relations: \vspace{-0.1cm}
$$\begin{array}{ccccc}
    zx=xz+y & \mbox{ } & yx=xy-2x & \mbox{ } & zy=yz-2z  \\
\end{array}$$
The next result can be found in $\cite{CK}$,$\cite{YS}$.
\begin{proposition}
 The map $\rho:\mathbb{C}\langle x,y,z \rangle /I_{{\mathfrak
{sl}}_2} \longrightarrow \en(\mathbb{C}[x_1,x_2])$ given by
$\rho(x)={\displaystyle  x_2 \frac{\partial}{ \partial x_1}}$,
$\rho(y)={\displaystyle  x_1 \frac{\partial}{ \partial x_1}-x_2
\frac{\partial}{ \partial x_2} } $ and
$\rho(z)={\displaystyle  x_1 \frac{\partial}{ \partial x_2}}$
defines a representation of the algebra $\mathbb{C}\langle x,y,z
\rangle /I_{{\mathfrak {sl}}_2}$.
 \end{proposition}

Given  $s,n\in \mathbb{N}$ with $0\leq s \leq n$, the $s$-th
elementary symmetric function  ${\des \sum_{1\leq
i_1<\dots<i_s\leq n} x_{i_1}\dots x_{i_s}}$ on variables
$x_1,\dots,x_n$ is denoted by $e_{s}^{n}(x_1,\dots,x_n)$. For
$b\in \mathbb{N}$, the notation
$e_{s}^{n}(b):=e_{s}^{n}(b,b-1,\dots, b-n+1)$  we will used. Given
$a,n\in \mathbb{N}$ such that $a\leq n$, we set $(a)_n=a(a-1)\dots
(a-n+1)$.

\begin{theorem}\label{ncsl2} Given $a,b\in \mathbb{N}$, the following
identities hold in $\mathbb{C}\langle x,y,z \rangle /I_{{\mathfrak
{sl}}_2}$
\begin{enumerate}
\item{$z^{a}x^{b}=\des \sum_{s,k}  \frac{(a)_k (b)_k}{k!} e_{k-s}^{k}(-a-b+2k)x^{b-k}
y^{s}z^{a-k}$, \\where the sum runs over all $k,s\in \mathbb{N}$
such that $0\leq s\leq k
\leq
\min(a,b)$.}
\item{$z^{a}y^{b}=\des{\sum_{k=0}^{b} {b\choose
k}(-2a)^{k}y^{b-k}z^{a} }$.}
\item{$y^{a}x^{b}=\des{\sum_{k=0}^{a} {a\choose
k}(-2b)^{k}x^{b}y^{a-k}}$.}
\end{enumerate}
\end{theorem}
\begin{proof}
Formula $1.$ is proved by induction. It is equivalent to another
formula for the normalization of $z^{a}x^{b}$ given in
\cite{Kac1}. $2.$ and $3.$ are similar to Theorem \ref{hwe}, part $1$. Formula
$2.$ express the fact that as we try to move the $z$'s above the
$y$'s some of the $y$'s may get killed. This argument justify the
${b\choose k}$ factor. The $a^{k}$ factor arises from the fact
that each `$y$' may be killed for any of the $z$'s. The $(-2)^{k}$
factor follows from the fact that each killing of a `$y$' is
weighted by a $-2$.
\end{proof}

Assume we are given $A=(A_1,\dots,A_n)\in (\mathbb{N}^{3})^{n}$
and $A_i=(a_i,b_i,c_i)$, for $i\in[[1,n]]$. Set
$X^{A_{i}}=x^{a_i}y^{b_i}z^{c_i}$ for $i\in[[1,n]]$, furthermore
set $a=(a_1,\dots,a_n)\in \mathbb{N}^{n}$, $b=(b_1,\dots,b_n)\in
\mathbb{N}^{n}$, $c=(c_1,\dots,c_n)\in \mathbb{N}^{n}$ and
$|A|=(|a|,|b|,|c|)\in \mathbb{N}^{3}$. Using this notation, we
have the\\
\begin{definition}
The {\em normal coordinates} $ N_{\mathfrak{sl}_2}(A,k,s,p,q)$ of
${\des\prod_{i=1}^{n}X^{A_{i}}}$ in the algebra ${\des
\mathbb{C}\langle x,y,z
\rangle /I_{{\mathfrak {sl}}_2}}$ are given via the identity
\begin{equation}\label{nc}
{\des\prod_{i=1}^{n}X^{A_{i}}=  \sum_{k,s,p,q}
N_{\mathfrak{sl}_2}(A,k,s,p,q) X^{r(A,k,s,p,q)}}
\end{equation}
where $k,s,p,q\in \mathbb{N}^{n-1}$ are such that $0\leq s_i\leq
k_i\leq \min(|c_{\leq i}|-|k_{<i}|,a_{i+1})$, \\ $0\leq p_i \leq
b_{i+1}$, $0\leq q_i\leq |b_{\leq i}|+|s_{<i}|-|p_{<i}|-|q_{<i}|$,
for  $i\in[[1,n-1]]$. Moreover,
$r(A,k,s,p,q)=(|a|-|k|,|b|+|s|-|p|-|q|,|c|-|k| ).$  We set
$N_{\mathfrak{sl}_2}(A,k,s,p,q)=0$ for $k,s,p,q$ not satisfying
the previous conditions.
\end{definition}
Theorem \ref{cnsl2} provides an explicit formula for the normal
coordinates $N_{\mathfrak{sl}_2}(A,k,s,t)$ of
${\des\prod_{i=1}^{n}X^{A_{i}}}$. Its proof goes by induction
using Theorem \ref{ncsl2}.

\begin{theorem}\label{cnsl2} With the notation of the definition above, we have
$$N_{\mathfrak{sl}_2}(A,k,s,p,q )= {\des (-2)^{|p|+|q|}
\prod_{i=1}^{n-1}\alpha_i \beta_i \gamma_i } {b_{i+1}\choose p_i} (|c_{\leq i}|-|k_{\leq
i}|)^{{p_i}}(a_{i+1}-k_i)^{q_i},
$$ where $\alpha_i={\des \frac{(|c_{\leq i}|- |k_{<i}|)_{k_i} (a_{i+1})_{k_i}}{k_i!}}$,
${\des\beta_i= e_{k_i-s_i}^{k_i}(-a_{i+1}- |c_{\leq i}|+|k_{<
i}|+2k_i)}$, \\ and ${\des\gamma_i={{|b_{\leq
i}|+|s_{<i}|-|p_{<i}|-|q_{<i}|}\choose q_{i}}}$, for all
$i\in[[1,n-1]]$.
\end{theorem}

Our final result provides an explicit formula for the product of
$m$ elements in the algebra $\sym^{n}(\mathbb{C}\langle x,y,z
\rangle /I_{{\mathfrak {sl}}_2})$. Fix  $A:[[1,m]]\times
[[1,n]]\longrightarrow \mathbb{N}^{3}$, with
$(A_{ij})=((a_{ij}),(b_{ij}),(c_{ij}))$. Recall that given
$\sigma\in \mathbb{S}_n^{m}$ and $j\in[[1,n]]$, $A_j^{\sigma}$
denotes the vector $(A_{1\sigma_1^{-1}(j)},\dots,
A_{m\sigma_m^{-1}(j)})\in (\mathbb{N}^{3})^{m}$. Set
$X_j^{A_{ij}}=x_j^{a_{ij}} y_j^{b_{ij}} z_j^{c_{ij}}$, for
$j\in[[1,n]]$ and $|A_j^{\sigma}|=(|a_j^{\sigma}|,
|b_j^{\sigma}|,|c_j^{\sigma}|)$, where ${\des
|a_j^{\sigma}|=\sum_{i=1}^{m} a_{i\sigma^{-1}_i(j)}}$ and
similarly for $|b_j^{\sigma}|$, $|c_j^{\sigma}|$, and $k,s,p,q\in
(\mathbb{N}^{m-1})^{n}$.

\begin{theorem}\label{fpsl2} For any $A:[[1,m]]\times
[[1,n]]\longrightarrow\mathbb{N}^{3}$, the identity
$${\des (n!)^{m-1}\prod_{i=1}^{m} \overline{ \prod_{j=1}^{n}
X_j^{A_{ij}}}= \sum_{\sigma,k,s,p,q}\left( \prod_{j=1}^{n}
N_{\mathfrak{sl}_2}(A_j^{\sigma},k^{j},s^{j},p^{j},q^{j})\right)
\overline{
\prod_{j=1}^{n} X_j^{r_j(A_j^{\sigma},k^{j},s^{j},p^{j},q^{j})}} }$$
holds in $\sym^{n}(\mathbb{C}\langle x,y,z \rangle /I_{{\mathfrak
{sl}}_2})$, where
$k=(k^{1},\dots,k^{n})\in(\mathbb{N}^{m-1})^{n}$, and similar for
$s,p,q$,
$r_j(A_j^{\sigma},k^{j},s^{j},p^{j},q^{j})=(|a_j^{\sigma}|-|k^{j}|,
|b_j^{\sigma}|+|s^{j}|-|p^{j}|-|q^{j}|,|c_{j}^{\sigma}|-|k^{j}|$,
and $\sigma\in
\{\id\}\times \mathbb{S}_n^{m-1}$.

\end{theorem}

Theorem  \ref{fpsl2} is proven similarly to Theorem \ref{jo}.

\subsection*{Acknowledgment}
We thank Nicolas Andruskiewitsch, Sylvie Paycha and Carolina
Teruel. We also thank the organizing committee of  Geometric and
Topological Methods for Quantum Field Theory, Summer School 2003,
Villa de Leyva, Colombia.

\bibliographystyle{plain}
\bibliography{sqwa}

\begin{thebibliography}{10}

\bibitem{Al}
{J}. {A}lev, {T}.{J}. {H}odges, and {J}.{D}. {V}elez.
\newblock Fixed rings of the {W}eyl algebra {$A_1(\mathbb{C})$}.
\newblock {\em Journal of algebra}, 130:83--96, 1990.

\bibitem{CK}
{C}hristian {K}assel.
\newblock {\em Quantum groups}.
\newblock Springer-Velarg, New York, 1995.

\bibitem{EM}
{E}mil {M}artinec and {G}regory {M}oore.
\newblock {N}oncommutative {S}olitons on {O}rbifolds.
\newblock hep-th/0101199, 2001.

\bibitem{Kon}
{M}. {K}ontsevich.
\newblock Deformation {Q}uantization of {P}oisson {M}anifolds {I}.
\newblock math. q-alg/97090401, 1997.

\bibitem{YS}
{L}eonid {K}orogodski and {Y}an {S}oibelman.
\newblock Algebras of functions on {Q}uantum groups. part {I}.
\newblock {\em Mathematical surveys and monographs}, 56, 1996.

\bibitem{Etg2}
{P}avel {E}tingof and {V}ictor {G}inzburg.
\newblock {S}ymplectic reflection algebras, {C}alogero-{M}oser space, and
  deformed {H}arish-{C}handra homomorphism.
\newblock {\em Invent. Math}, 147(2):243--348, 2002.

\bibitem{GR}
{P}eter {D}oubilet.
\newblock {O}n the foundations of combinatorial theory. {VII}: {S}ymmetric
  functions through the theory of distribution and occupancy.
\newblock {\em Gian-Carlo Rota on Combinatorics. Introductory papers and
  commentaries}, pages 402--422, 1995.

\bibitem{DP}
{R}afael {D}\'{\i}az and {E}ddy {P}ariguan.
\newblock Quantum symmetric functions.
\newblock math.QA/0312494, 2003.

\bibitem{DP1}
{R}afael {D}\'{\i}az and {E}ddy {P}ariguan.
\newblock Super, quantum and non-commutative species.
\newblock Work in progress, 2004.

\bibitem{Min2}
{R}ajesh {G}opakumar, {S}hiraz {M}inwalla, and {A}ndrew {S}trominger.
\newblock {N}oncommutative {S}olitons.
\newblock {\em J. High Energy Phys. JHEP}, 05-020, 2000.

\bibitem{Sol}
{A}.{I}. {S}olomon.
\newblock {\em {P}hys.{L}ett. {A} 196}, 126(29), 1994.

\bibitem{Kac1}
{V}ictor {K}ac.
\newblock {\em Infinite dimensional Lie algebras.}
\newblock {C}ambridge {U}niversity {P}ress, New York, 1990.

\bibitem{Kac}
{V}ictor {K}ac and {P}okman {C}heung.
\newblock {\em {Q}uantum {C}alculus}.
\newblock Springer-Velarg, New York, 2002.

\end{thebibliography}

$$\begin{array}{c}
\mbox{Rafael D\'\i az. Instituto Venezolano de Investigaciones Cient\'\i ficas.} \  \mbox{\texttt{radiaz@ivic.ve}} \\
\!\!\!\!\!\mbox{Eddy Pariguan. Universidad central de Venezuela.} \  \mbox{\texttt{eddyp@euler.ciens.ucv.ve}} \\
\end{array}$$

\end{document}